\documentclass[12 pt]{amsart}
\usepackage{amscd,amsfonts,amssymb,amsmath}
\newtheorem{theorem}{Theorem}[section]
\theoremstyle{definition}

\newtheorem{lemma}[theorem]{Lemma}
\newtheorem{proposition}[theorem]{Proposition}

\newtheorem{problem}[theorem]{Problem}
\newtheorem{definition}[theorem]{Definition}
\newtheorem{example}[theorem]{Example}
\newtheorem{remark}[theorem]{Remark}
\numberwithin{equation}{subsection}
\usepackage{graphicx}

\newcommand{\Aut}{\operatorname{Aut}}
\newcommand{\End}{\operatorname{End}}

\begin{document}
\title[Singular braids]{Singular braids, singular links and subgroups of camomile type}
\author{Valeriy G. Bardakov}
\author{Tatyana A. Kozlovskaya}

\address{Sobolev Institute of Mathematics, 4 Acad. Koptyug avenue, 630090, Novosibirsk, Russia.}
\address{Novosibirsk State Agrarian University, Dobrolyubova street, 160, Novosibirsk, 630039, Russia.}
\address{Regional Scientific and Educational Mathematical Center of Tomsk State University, 36 Lenin Ave., Tomsk, Russia.}
\email{bardakov@math.nsc.ru}

\address{Regional Scientific and Educational Mathematical Center of Tomsk State University, 36 Lenin Ave., Tomsk, Russia.}
\email{t.kozlovskaya@math.tsu.ru}


\begin{abstract}
In this paper we find a finite set of generators and defining relations for the singular pure braid group $SP_n$, $n \geq  3$, that is a subgroup of the singular braid group $SG_n$.
Using this presentation, we prove that the center of $SG_n$ (which is equal to the center of $SP_n$ for  $n \geq 3$) is a direct factor in $SP_n$ but it is not a direct factor in $SP_n$. 
We introduce subgroups of camomile type and prove that the singular pure braid group $SP_n$, $n \geq  5$, is a subgroup of camomile type in $SG_n$.
Also we  construct the fundamental singquandle using a  representation of the singular braid monoid by endomorphisms of free guandle. For any singular link
we define some family of groups which are invariants of this link.
\end{abstract}

\maketitle

\section{Introduction}

Singular knot theory have gained a lot of interest in the past  decades. This was primarily motivated by Vassiliev invariants (finite  type invariants) \cite{V}.  Under the influence of the theory of Vassiliev invariants   singular braids were introduced. The relation between singular knots and singular braids is just the same as in the classical case. It is a natural problem to study their algebraic and geometric properties to construct invariants of singular links.

The singular braids with $n$ strands were introduced independently by J. Baez in \cite{Baez} and J. Birman in  \cite{Bir}.
These singular braids form a monoid $SB_n$ that is generated by the standard generators $\sigma_1^{\pm 1} ,\ldots, \sigma_{n-1}^{\pm 1} $ of the braid group $B_n$ plus the additional singular generators $\tau_1 ,\ldots, \tau_{n-1}$.
 It is shown in \cite{FKR} that  the Baez-Birman monoid on $n$ strands $SB_n$ embeds into a group $SG_n$ that is now known as the singular braid group on $n$ strands.
The word problem for $SB_n$ and $SG_n$ was solved in  \cite{Co, Par, V4}.
We refer to  \cite{DG, DG1, V3} for more on singular braid monoid and singular braid groups.

The pure braid group $P_n$ is the kernel of the epimorphism  $B_n$ to the symmetric group $S_n$. The singular pure braid groups $SP_n$  are generalizations of $P_n$.  $SP_n$ was introduced  in   \cite{DG}.  The group $SP_n$ is the kernel of the epimorphism that maps, for each $i$,  $\sigma_i$ and $\tau_i$ to the cyclic permutation $(i, i+1)$. Dasbach and Gemein  \cite{DG} found a set of generators and defining relations for this group and established that this group can be constructed using successive HNN-extensions.  A decomposes $SP_3$ and $SP_4$ as a semi-direct product of two groups was revisited in \cite{BK} and \cite{K-4}.
In \cite{DG1} Dasbach and Gemein investigated extensions of the Artin representation $B_n \to \Aut(F_n)$ and the Burau representation $B_n \to GL_n(\mathbb{Z}[t, t^{-1}])$ to the singular braid monoid and found connections between these representations.  They showed  that a certain linear representation of the singular braid monoid  $SB_3$ is faithful.

Also, it is possible to define other homomorphism of $SG_n$ to $S_n$, which sends $\sigma_i$ to $(i, i+1)$ and $\tau_i$ to $e$. The kernel of this homomorphism is denoted by $ST_n$.
In \cite{GKM}, it was found a presentation for the group $ST_3$ and proved  that it is isomorphic to $SP_3$.

A lot of papers are dedicated to construction of invariants for singular links.
The HOMFLY and Kauffman polynomials were extended to 3-variable polynomials for
singular links by Kauffman and Vogel \cite{KV}. The extended HOMFLY polynomial was recovered by the construction of traces on singular Hecke
algebras \cite{PR}.
Juyumaya and Lambropoulou  \cite{JL} used a similar approach to define invariants of singular links.

 The Khovanov homology was extended to a homology for singular links in \cite{Sh}. The Alexander polynomials of a cube of resolutions (in Vassiliev's sense) of a singular
knot were categorified in \cite{A}. Moreover, a 1-variable extension of the Alexander polynomial for singular links was categorified in \cite{OSS}. The generalized
cube of resolutions (containing Vassilievs resolutions as well as those smoothings at double points which preserve the orientation) was categorified in \cite{OS}.

Fiedler \cite{Fd} extended the Kauffman state models of the Jones and Alexander polynomials to the context of singular knots. 
He extends the Kauffman state models of the Jones and Alexander
polynomials of classical links to state models of their two-variable
extensions in the case of singular links.

The theory of singular braid is connected with the theory of pseudo-braids. In particular, 
in  \cite{BJW} was proved that the monoid of pseudo-braids is isomorphic to the singular braid monoid.
Hence, the group of the singular braids is isomorphic to the group of pseudo-braids. On the other side, the theory of pseudo-links is a quotient of the theory of singular  links by the singular first Reidemeister move.  In \cite{NOS} was generalize the notion of quandles to psyquandles and used these to define invariants of oriented
singular links and pseudo-links. In particular, in this paper was introduced Alexander
psyquandles and a generalization of the Alexander polynomial for oriented singular links and pseudo-links.

In this paper we study the structure of the  singular pure braid group $SP_n$.
We find a finite set of generators and defining relations for  $SP_n$, $n \geq  3$.
It is known that the center of $P_n$ is equal to the center of $B_n$, $n>2$. Vershinin \cite{V4} proved that $Z(B_n) = Z(SP_n) = Z(SG_n)$ is infinite cyclic group.
Using the presentation of $SP_n$, we prove that for $n > 2$ the center $Z(SP_n)$  is a direct factor in $SP_n$, but it not not a direct factor in $SG_n$. The similar result for the center of $P_n$ and $B_n$ was proved by Neshchadim \cite{N1, N2}.

During last decades people study groups which look like  the braid group $B_n$ (see \cite{B,  Bir,  KL, V3}).
For example, virtual braid group $VB_n$, welded  braid group $WB_n$,
singular  braid group $SG_n$ and some other.
Any group of this type  has a pure subgroup, which is the kernel of an epimorphism onto $S_n$. The structures of these pure subgroups are different for the case $n=3$ and $n>3$. In the second case we have far commutativity relations which we do not have in the first case. The set of generators and defining relations any pure group for  $n>4$ can be gotten  from the set of generators and defining relations of the pure group for  $n=4$, using conjugations by elements of the corresponding braid group. Hence, in some sense it is need to study the case $n=3$ and the case $n=4$.  Formalizing this observation,
we introduce subgroups of camomile type and prove that the singular pure braid group $SP_n$, $n \geq  5$, is a subgroup of camomile type in $SG_n$.

In the  classical knot theory for any diagram of a link can be defined a quandle which is a link invariant (see \cite{M}, \cite{J}). 
In \cite{E} was introduced  an  oriented singquandle whose axioms  come  from singular  Reidemeister moves on oriented singular link diagram. Using a  representation of the singular braid monoid by endomorphisms of free singguandle, we suggest another approach for 
construction fundamental singquandles for singular links. For any singular link we define some family of groups which are invariants of this link.

The paper is organizes as follows. In Section \ref{BD} wee recall some known facts from the theory of classical braids and singular braids. 
In Section \ref{pure} we are studying the singular pure braid group  $SP_n$. We give two Schreier sets of coset representatives of $SP_n$ in $SG_n$,
find  a finite presentation for $SP_n$, $n \geq 2$, which is simpler than the representation in  \cite{DG}. At the end of this section we give conjugation 
rules of actions by generators $SG_n$ on the generators  $SP_n$.
In Section \ref{cam} we give a definition subgroup of camomile type and prove that any $SP_n$ for $n\geq 5$ is a subgroup of camomile type in $SG_n$ with highlighted petal $SP_4$.
In Section \ref{cent} we prove that the center of $SP_n$, $n \geq 3$ (which is infinite cyclic and is equal to the center of $SB_n$) is a direct factor. On the other side 
we show that it is not a direct factor in $SG_n$.
In Section \ref{Gr} we  consider
singquandles which was introduced in \cite{E} and studied in \cite{BC, E1}. For any singular link $L$ in \cite{E} was define a fundamental oriented singquandle  and was prove that it is an invariant.
We give a representation of $SB_n$ by endomorphisms of free singquandle $FSQ_n$ and using this representation we   give another interpretation of fundamental oriented singquandle.
In Section \ref{GSL},  using representations of $SB_n$ by endomorphisms of free group, which was constructed in \cite{E}, we construct some  group invariants  of singular links.

In the last section we formulate some open problems and suggest directions for further research.

\section*{Acknowledgments}
This work  is supported by the Ministry of Science and Higher Education of Russia (agreement  No.  075-02-2022-884)

\bigskip

\section{Basic definitions} \label{BD}

In this section we recall some known definitions which can be found in \cite{Artin, Bir1, Mar}. 

\begin{figure}[h]
\includegraphics[totalheight=8cm]{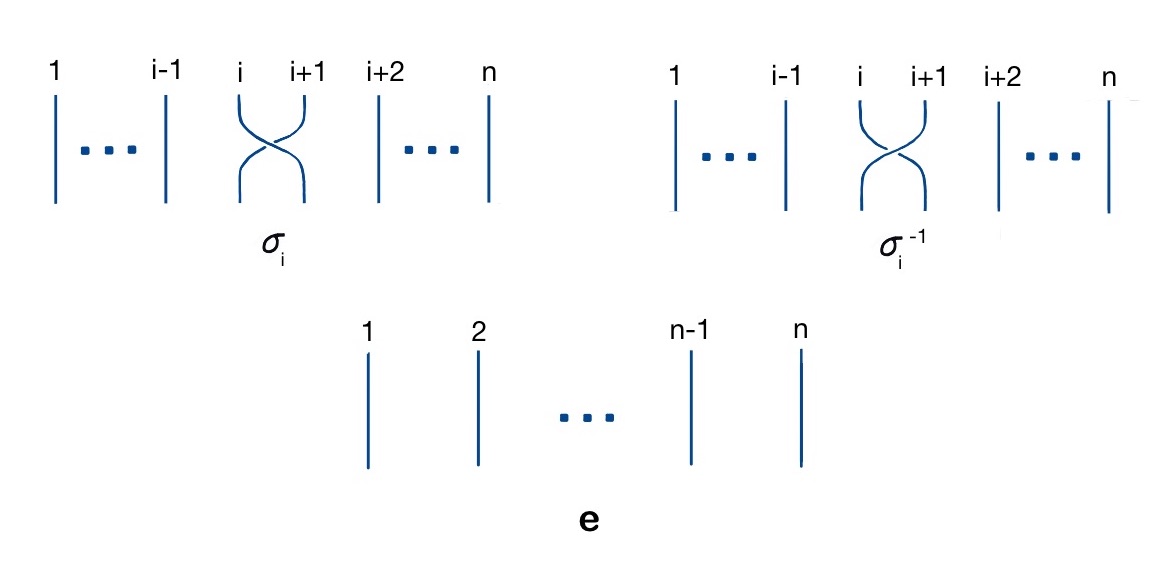}
\caption{The elementary braids $\sigma_i$, $\sigma_i^{-1}$ and the unit $e$} \label{figure1}
\end{figure}

The braid group $B_n$, $n\geq 2$, on $n$ strands can be defined as
a group generated by $\sigma_1,\sigma_2,\ldots,\sigma_{n-1}$ with the defining relations
\begin{equation}
\sigma_i \, \sigma_{i+1} \, \sigma_i = \sigma_{i+1} \, \sigma_i \, \sigma_{i+1},~~~ i=1,2,\ldots,n-2, \label{eq1}
\end{equation}
\begin{equation}
\sigma_i \, \sigma_j = \sigma_j \, \sigma_i,~~~|i-j|\geq 2. \label{eq2}
\end{equation}
The geometric interpretation of  $\sigma_i$, its inverse $\sigma_{i}^{-1}$ and the unit $e$ of $B_n$ are depicted  in the Figure~\ref{figure1}.

There exists a homomorphism of $B_n$ onto the symmetric group $S_n$ on
$n$ symbols. This homomorphism  maps
 $\sigma_i$ to the transposition  $(i,i+1)$, $i=1,2,\ldots,n-1$.
The kernel of this homomorphism is called the
{\it pure braid group} and denoted by
$P_n$. The group $P_n$ is generated by  $a_{ij}$, $1\leq i < j\leq n$.
These
generators can be expressed by the generators of
 $B_n$ as follows
$$
a_{i,i+1}=\sigma_i^2,
$$
$$
a_{ij} = \sigma_{j-1} \, \sigma_{j-2} \ldots \sigma_{i+1} \, \sigma_i^2 \, \sigma_{i+1}^{-1} \ldots
\sigma_{j-2}^{-1} \, \sigma_{j-1}^{-1},~~~i+1< j \leq n.
$$
In these generators $P_n$ is defined by relations
 \begin{align}
& a_{ik} a_{ij} a_{kj} = a_{kj} a_{ik} a_{ij},   \label{re2}\\
& a_{mj} a_{km} a_{kj} = a_{kj} a_{mj} a_{km},  ~\mbox{for}~m < j, \label{re3}\\
& (a_{km} a_{kj} a_{km}^{-1}) a_{im} = a_{im} (a_{km} a_{kj} a_{km}^{-1}),  ~\mbox{for}~i < k < m < j, \label{re4}\\
& a_{kj} a_{im} = a_{im} a_{kj},  ~\mbox{for}~k < i < m < j ~\mbox{or}~m < k. \label{re1}
\end{align}
The subgroup $P_n$ is normal in $B_n$, and the quotient $B_n / P_n$ is  $S_n$. The generators of $B_n$ act on the generator $a_{ij} \in P_n$ by the rules:
 \begin{align}
& \sigma_k^{-1} a_{ij} \sigma_k =  a_{ij},  ~\mbox{for}~k \not= i-1, i, j-1, j, \label{c1}\\
& \sigma_{i}^{-1} a_{i,i+1} \sigma_{i} =  a_{i,i+1},   \label{c2}\\
& \sigma_{i-1}^{-1} a_{ij} \sigma_{i-1} =   a_{i-1,j},   \label{c3}\\
& \sigma_{i}^{-1} a_{ij} \sigma_{i} =  a_{i+1,j} [a_{i,i+1}^{-1}, a_{ij}^{-1}],  ~\mbox{for}~j \not= i+1 \label{c4}\\
& \sigma_{j-1}^{-1} a_{ij} \sigma_{j-1} =  a_{i,j-1},   \label{c5}\\
& \sigma_{j}^{-1} a_{ij} \sigma_{j} =  a_{ij} a_{i,j+1} a_{ij}^{-1},   \label{c6}
\end{align}
where $[a, b] = a^{-1} b^{-1} a b = a^{-1} a^b$.

Denote by
$$
U_{i} = \langle a_{1i}, a_{2i}, \ldots, a_{i-1,i} \rangle,~~~i = 2, \ldots, n,
$$
a subgroup of $P_n$.
It is known that $U_i$ is a free group of rank $i-1$. One can rewrite the defining relations of $P_n$ as the following conjugation rules (for $\varepsilon = \pm 1$):
  \begin{align}
& a_{ik}^{-\varepsilon} a_{kj}  a_{ik}^{\varepsilon} = (a_{ij} a_{kj})^{\varepsilon} a_{kj} (a_{ij} a_{kj})^{-\varepsilon},  \label{co1}\\
& a_{km}^{-\varepsilon} a_{kj}  a_{km}^{\varepsilon} = (a_{kj} a_{mj})^{\varepsilon} a_{kj} (a_{kj} a_{mj})^{-\varepsilon},  ~\mbox{for}~m < j, \label{co2}\\
& a_{im}^{-\varepsilon} a_{kj}  a_{im}^{\varepsilon} = [a_{ij}^{-\varepsilon}, a_{mj}^{-\varepsilon}]^{\varepsilon} a_{kj} [a_{ij}^{-\varepsilon}, a_{mj}^{-\varepsilon}]^{-\varepsilon},  ~\mbox{for}~i < k < m, \label{co3}\\
& a_{im}^{-\varepsilon} a_{kj} a_{im}^{\varepsilon} = a_{kj},  ~\mbox{for}~k < i < m < j ~\mbox{or}~  m < k. \label{co4}
\end{align}
The group $P_n$ is a semi--direct product of  the normal subgroup
$U_n$ and the group $P_{n-1}$. Similarly, $P_{n-1}$ is a semi--direct product of the free group
$U_{n-1}$  and the group $P_{n-2},$ and so on.
Therefore, $P_n$ is decomposable (see \cite{Mar}) into the following semi--direct product
$$
P_n=U_n\rtimes (U_{n-1}\rtimes (\ldots \rtimes
(U_3\rtimes U_2))\ldots),~~~U_i\simeq F_{i-1}, ~~~i=2,3,\ldots,n.
$$

The group $B_n$ has a faithful representation into the automorphism group  ${\rm
Aut}(F_n)$ of the free group $F_n = \langle x_1, x_2, \ldots, x_n \rangle.$
In this case the generator $\sigma_i$, $i=1,2,\ldots,n-1$, goes to the automorphism
$$
\sigma_{i} : \left\{
\begin{array}{ll}
x_{i} \longmapsto x_{i} \, x_{i+1} \, x_i^{-1}, &  \\ x_{i+1} \longmapsto
x_{i}, & \\ x_{l} \longmapsto x_{l}, &  l\neq i,i+1.
\end{array} \right.
$$

By theorem of Artin \cite[Theorem 1.9]{Bir}, an automorphism $\beta $ from ${\rm
Aut}(F_n)$ lies in  $B_n$ if and only if $\beta $ satisfies  the following conditions:
$$
~~~~~1)~~ \beta(x_i) = a_i^{-1} \, x_{\pi(i)} \, a_i,~~~1\leq i\leq n,
$$
$$
2)~~ \beta(x_1x_2 \ldots x_n)=x_1x_2 \ldots x_n,
$$
where $\pi $ is a permutation from $S_n$ and $a_i\in F_n$.

{\it The Baez--Birman monoid}
\cite{Baez, Bir} or {\it the singular braid monoid} $SB_n$ is generated
(as a monoid) by elements $\sigma_i,$ $\sigma_i^{-1}$, $\tau_i$, $i = 1, 2, \ldots, n-1$.
The elements $\sigma_i,$ $\sigma_i^{-1}$  generate the braid group
$B_n$. The generators $\tau_i$  satisfy the defining relations
\begin{equation}
\tau_i \, \tau_j = \tau_j \, \tau_i, ~~~|i - j| \geq 2, \label{eq12}
\end{equation}
other relations are mixed:
\begin{equation}
\tau_{i} \, \sigma_{j} = \sigma_{j} \, \tau_{i}, ~~~|i - j| \geq 2, \label{eq13}
\end{equation}
\begin{equation}
\tau_{i}  \, \sigma_{i} = \sigma_{i} \, \tau_{i},~~~ i=1,2,\ldots,n-1,  \label{eq14}
\end{equation}
\begin{equation}
\sigma_{i} \, \sigma_{i+1} \, \tau_i = \tau_{i+1} \, \sigma_{i}  \, \sigma_{i+1},~~~ i=1,2,\ldots,n-2,
 \label{eq15}
 \end{equation}
 \begin{equation}
\sigma_{i+1}  \, \sigma_{i} \, \tau_{i+1} = \tau_{i} \,
\sigma_{i+1} \, \sigma_{i}, ~~~ i=1,2,\ldots,n-2.
 \label{eq16}
\end{equation}

\begin{figure}[h]
\includegraphics[width=7cm, height=3.6cm]{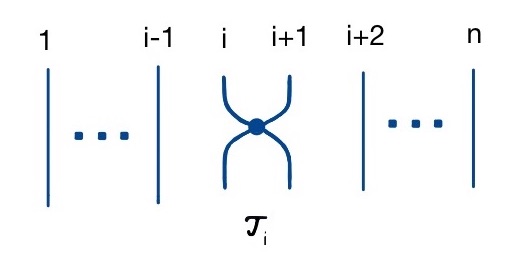}\\

\caption{The elementary singular braid  $\tau_i$} \label{figure}

\end{figure}

For a geometric interpretation of the elementary singular braid $\tau_i$ see Figure~\ref{figure}.

It is proved by R. Fenn, E. Keyman and C. Rourke \cite{FKR} that the Baez-Birman monoid $SB_n$  is embedded into a group $SG_n$  which they call the  {\it singular braid group}.

\bigskip

\section{Singular pure braid group}\label{pure}

Define the map
$$
\pi : SG_n \longrightarrow S_n
$$
of $SG_n$ onto the symmetric group $S_n$ on $n$ symbols by actions on the generators
$$
\pi(\sigma_i) = \pi(\tau_i) = (i, i+1), ~~~ i = 1, 2, \ldots, n-1.
$$
The kernel $\mbox{ker}(\pi)$ of this map is called the
{\it singular pure braid group} and is denoted by $SP_n$ (see \cite{DG}).
It is clear that $SP_n$ is a normal subgroup of index $n!$ of $SG_n$ and we have a short exact sequence
$$
1 \to SP_n \to SG_n \to S_n \to 1.
$$

Let $m_{kl} = \sigma_{k-1}  \, \sigma_{k-2} \ldots \sigma_l$ for $l < k$ and $m_{kl} = 1$
in other cases. Then the set
$$
\Lambda_n = \left\{ \prod\limits_{k=2}^n m_{k,j_k} \vert 1 \leq j_k
\leq k \right\}
$$
is a Schreier set of coset representatives of $SP_n$ in $SG_n$.

Also, we will use another set of coset representatives of $SP_n$ in $SG_n$. Put  $n_{kl} = \sigma_{k-1}^{-1}  \, \sigma_{k-2}^{-1}  \ldots \sigma_l^{-1}$ for $l < k$ and $n_{kl} = 1$
in other cases. Then the set
$$
M_n = \left\{ \prod\limits_{k=2}^n n_{k,j_k} \vert 1 \leq j_k
\leq k \right\}
$$
is a Schreier set of coset representatives of $SP_n$ in $SG_n$.

Define the following subsets of $M_n$, $n > 2$,
$$
M_{n,m} = M_n \setminus \left( M_m \setminus \{ e \} \right) ~\mbox{for}~1< m< n.
$$
In particular,
$$
M_{n,n-1} = M_n \setminus \left( M_{n-1} \setminus \{ e \} \right) = \{ e, \,  \sigma_{n-1}^{-1},  \, \sigma_{n-1}^{-1} \sigma_{n-2}^{-1},  \,  \ldots,  \, \sigma_{n-1}^{-1} \sigma_{n-2}^{-1} \ldots \sigma_1^{-1} \}.
$$

\subsection{Presentation of $SP_n$ }\label{gen}

 Dasbach and Gemein \cite{DG} have found some set of generators and defining relation for $SP_n$. We simplify this presentation and find some concrete finite representation (for $n=3,4$ it was done  in \cite{BK}, \cite{K-4}).
To describe it let us define the following elements in  $SP_n$:
$$
a_{i,i+1}=\sigma_i^2,~~~b_{i,i+1} = \tau_i \sigma_i,
$$
$$
a_{ij} = \sigma_{j-1} \, \sigma_{j-2} \ldots \sigma_{i+1} \, a_{i,i+1} \, \sigma_{i+1}^{-1} \ldots
\sigma_{j-2}^{-1} \, \sigma_{j-1}^{-1},~~~i+1< j \leq n,
$$
$$
b_{ij} = \sigma_{j-1} \, \sigma_{j-2} \ldots \sigma_{i+1} \, b_{i,i+1} \, \sigma_{i+1}^{-1} \ldots
\sigma_{j-2}^{-1} \, \sigma_{j-1}^{-1},~~~i+1< j \leq n.
$$

Using this notations we can prove the main result of the present section.

\begin{theorem} \label{mt}
The singular pure braid group $SP_n$, $n \ge 2$ is generated by elements ${a_{ij}}$, ${b_{ij}}$ %
$1 \le i < j \le n$ and is defined by relations ($\varepsilon=\pm1$):

\begin{equation} \label{dr-0}
 a_{ij} b_{ij} = b_{ij} a_{ij},
\end{equation}
\begin{equation} \label{dr-1}
a_{ik}^{-\varepsilon} a_{kj}  a_{ik}^{\varepsilon} = (a_{ij} a_{kj})^{\varepsilon} a_{kj} (a_{ij} a_{kj})^{-\varepsilon},
\end{equation}
\begin{equation} \label{dr-11}
a_{ik}^{-\varepsilon} b_{kj}  a_{ik}^{\varepsilon} = (a_{ij} a_{kj})^{\varepsilon} b_{kj} (a_{ij} a_{kj})^{-\varepsilon},
\end{equation}
\begin{equation} \label{dr-2}
 a_{km}^{-\varepsilon} a_{kj}  a_{km}^{\varepsilon} = (a_{kj} a_{mj})^{\varepsilon} a_{kj} (a_{kj} a_{mj})^{-\varepsilon},  ~\mbox{for}~m < j,
\end{equation}
\begin{equation} \label{dr-21}
 a_{km}^{-\varepsilon} b_{kj}  a_{km}^{\varepsilon} = (a_{kj} a_{mj})^{\varepsilon} b_{kj} (a_{kj} a_{mj})^{-\varepsilon},  ~\mbox{for}~m < j,
\end{equation}
\begin{equation} \label{dr-3}
a_{im}^{-\varepsilon} a_{kj}  a_{im}^{\varepsilon} = [a_{ij}^{-\varepsilon}, a_{mj}^{-\varepsilon}]^{\varepsilon} a_{kj} [a_{ij}^{-\varepsilon}, a_{mj}^{-\varepsilon}]^{-\varepsilon},  ~\mbox{for}~i < k < m,
\end{equation}
\begin{equation} \label{dr-31}
a_{im}^{-\varepsilon} b_{kj}  a_{im}^{\varepsilon} = [a_{ij}^{-\varepsilon}, a_{mj}^{-\varepsilon}]^{\varepsilon} b_{kj} [a_{ij}^{-\varepsilon}, a_{mj}^{-\varepsilon}]^{-\varepsilon},  ~\mbox{for}~i < k < m,
\end{equation}
\begin{equation} \label{dr-4}
 a_{im}^{-\varepsilon} a_{kj} a_{im}^{\varepsilon} = a_{kj},  ~\mbox{for}~k < i < m < j ~\mbox{or}~  m < k,
\end{equation}
\begin{equation} \label{dr-41}
 a_{im}^{-\varepsilon} b_{kj} a_{im}^{\varepsilon} = b_{kj},  ~\mbox{for}~k < i < m < j ~\mbox{or}~  m < k,
\end{equation}
\begin{equation} \label{dr-5}
 b_{im}^{-\varepsilon} a_{kj} b_{im}^{\varepsilon} = a_{kj},  ~\mbox{for}~k < i < m < j ~\mbox{or}~  m < k,
\end{equation}
\begin{equation} \label{dr-51}
 b_{im}^{-\varepsilon} b_{kj} b_{im}^{\varepsilon} = b_{kj},  ~\mbox{for}~k < i < m < j ~\mbox{or}~  m < k,
\end{equation}
\begin{equation} \label{dr-6}
b_{ij}^{-\varepsilon} (a_{ik} a_{jk}) b_{ij}^{\varepsilon} = a_{ik} a_{jk}, ~\mbox{for}~i < j < k,
\end{equation}
\begin{equation} \label{plast}
b_{im}^{-\varepsilon} \left(a_{mj}^{-1} a_{kj}  a_{mj} \right)  b_{im}^{\varepsilon}= a_{mj}^{-1} a_{kj}  a_{mj},  ~\mbox{for}~i < k < m,
\end{equation}
\begin{equation} \label{last}
b_{im}^{-\varepsilon} \left(a_{mj}^{-1} b_{kj}  a_{mj} \right)  b_{im}^{\varepsilon}= a_{mj}^{-1} b_{kj}  a_{mj},  ~\mbox{for}~i < k < m.
\end{equation}
\end{theorem}

\begin{proof}
We use the standard Reidemester-Shraier method (see for example \cite[Section~2.3]{MKS}). If we consider a subgroup of $SP_n$, that is generated by $a_{ij}$, we get the pure braid group $P_n$. We know defining relations of this group   (see Section \ref{BD}). The generators and relations  of $SP_3$ have found in \cite{BK}. By analogy we can find relations of $SP_n$ in general case.

We discuss  only two last type of relations, which different from relations of $P_n$ and do not arrive in the presentation of $SP_3$.

From the commutativity relation
$$
\sigma_j^{-1} \tau_i^{-1} \sigma_j \tau_i = 1,~~|j-i| >1,
$$
which hold in $SG_n$, $n \geq 4$, follows the commutativity relation
$$
[a_{j,j+1}, b_{i,i+1}] = 1
$$
in  $SP_n$, $n \geq 4$. Suppose that $j > i + 1$ and conjugating this relation by $\sigma_{j-1}^{-1} \sigma_{j-2}^{-1} \ldots \sigma_{i+1}^{-1}$, we get
$$
[a_{j,j+1}, b_{i,i+1}]^{\sigma_{j-1}^{-1} \sigma_{j-2}^{-1} \ldots \sigma_{i+1}^{-1}} = [a_{i+1,j+1}^{a_{i+2,j+1}(a_{i+3,j+1} \ldots a_{j,j+1})}, b_{i,i+2}].
$$
Since
$$
\left( a_{i+3,j+1} \ldots a_{j,j+1} \right)^{b_{i,i+2}} = a_{i+3,j+1} \ldots a_{j,j+1},
$$
we have
\begin{equation} \label{re-1}
[a_{i+1,j+1}^{a_{i+2,j+1}}, b_{i,i+2}] = 1.
\end{equation}
This relation can be write in the form
$$
(a_{i+2,j+1}^{-1} a_{i+1,j+1} a_{i+2,j+1})^{b_{i,i+2}} = a_{i+2,j+1}^{-1} a_{i+1,j+1} a_{i+2,j+1}.
$$
It is relation of the type (\ref{plast}).

If $j > i+2$, then we take $k$ such that $i+2 \leq k \leq j-1$. Conjugating the relation (\ref{re-1}) by $\sigma_{i+2}^{-1} \sigma_{i+3}^{-1} \ldots \sigma_{k}^{-1}$, we get
$$
[a_{i+1,j+1}^{a_{i+2,j+1}}, b_{i,i+2}]^{\sigma_{i+2}^{-1} \sigma_{i+3}^{-1} \ldots \sigma_{k}^{-1}} = [a_{i+1,j+1}^{a_{k+1,j+1}}, b_{i,k+1}] = 1.
$$
That is equivalent to the relation
\begin{equation} \label{re-2}
(a_{k+1,j+1}^{-1} a_{i+1,j+1} a_{k+1,j+1})^{b_{i,k+1}} = a_{k+1,j+1}^{-1} a_{i+1,j+1} a_{k+1,j+1}.
\end{equation}
It is relation of the type (\ref{plast}). In particular, for $k= j-1$, we have
$$
[a_{i+1,j+1}^{a_{j,j+1}}, b_{i,j}] = 1.
$$
If we conjugate it by $\sigma_j^{-1}$, we have
$$
[a_{i+1,j+1}^{a_{j,j+1}}, b_{i,j}] ^{\sigma_{j}^{-1}} = [a_{i+1,j}, b_{i,j+1}] = 1.
$$
That is equivalent to the relation
$$
b_{i,j+1}^{a_{i+1,j}} = b_{i,j+1}
$$
of the form (\ref{dr-5}).

Conjugating the defining relations (\ref{re-2}) by elements $\lambda^{-1}$, $\lambda \in \Lambda_n$, we get other relations of the form (\ref{plast}).

To get relations (\ref{last}) we take the commutativity relation
$$
\tau_j^{-1} \tau_i^{-1} \tau_j \tau_i = 1,~~|j-i| >1,
$$
which hold in $SG_n$, $n \geq 4$. From this relation  follows the commutativity relation
$$
[b_{j,j+1}, b_{i,i+1}] = 1
$$
and we can construct relations (\ref{last}), using the same approach as before.
\end{proof}

Subgroup $SP_n$ is normal in $SG_n$ and one can find the formulas of conjugations of the generators of $SP_n$ by the generators  $\sigma_i$, $i=1, 2, \ldots, n-1$ (for $\varepsilon = \pm 1$),

 \begin{align}
& \sigma_k^{-\varepsilon} a_{ij} \sigma_k^{\varepsilon} =  a_{ij}, ~~ \sigma_k^{-\varepsilon} b_{ij} \sigma_k^{\varepsilon} =  b_{ij}~~~\mbox{for}~k \not= i-1, i, j-1, j, \label{ct1}\\
& \sigma_{i}^{-\varepsilon} a_{i,i+1} \sigma_{i}^{\varepsilon} =  a_{i,i+1},~~  \sigma_{i}^{-\varepsilon} b_{i,i+1} \sigma_{i}^{\varepsilon} =  b_{i,i+1}, \label{ct2}\\
& \sigma_{i-1}^{-1} a_{ij} \sigma_{i-1} =   a_{i-1,j},~~  \sigma_{i-1}^{-1} b_{ij} \sigma_{i-1} =   b_{i-1,j},  \label{ct3}\\
& \sigma_{i-1} a_{ij} \sigma_{i-1}^{-1} =   a_{ij}^{-1} a_{i-1,j} a_{ij},~~  \sigma_{i-1} b_{ij} \sigma_{i-1}^{-1} =  a_{ij}^{-1}  b_{i-1,j} a_{ij},  \label{ct31}\\
& \sigma_{i}^{-1} a_{ij} \sigma_{i} =  a_{ij} a_{i+1,j} a_{ij}^{-1}, ~~ \sigma_{i}^{-1} b_{ij} \sigma_{i} =  a_{ij} b_{i+1,j} a_{ij}^{-1} ~~~\mbox{for}~j \not= i+1 \label{ct4}\\
& \sigma_{i} a_{ij} \sigma_{i}^{-1} =   a_{i+1,j},~~  \sigma_{i} b_{ij} \sigma_{i}^{-1} = b_{i+1,j} ~~~\mbox{for}~j \not= i+1 \label{ct41}\\
& \sigma_{j-1}^{-1} a_{ij} \sigma_{j-1} =  a_{i,j-1}, ~~   \sigma_{j-1}^{-1} b_{ij} \sigma_{j-1} =  b_{i,j-1}, \label{ct5}\\
& \sigma_{j-1} a_{ij} \sigma_{j-1}^{-1} =   a_{j-1,j} a_{i,j-1} a_{j-1,j}^{-1},~~  \sigma_{j-1} b_{ij} \sigma_{j-1}^{-1} =   a_{j-1,j} b_{i,j-1} a_{j-1,j}^{-1}, \label{ct51}\\
& \sigma_{j}^{-1} a_{ij} \sigma_{j} =  a_{j,j+1}^{-1} a_{i,j+1} a_{j,j+1}, ~~  \sigma_{j}^{-1} b_{ij} \sigma_{j} =  a_{j,j+1}^{-1} b_{i,j+1} a_{j,j+1}
~~~\mbox{for}~1 \leq i < j \leq n-1, \label{ct6}\\
& \sigma_{j} a_{ij} \sigma_{j}^{-1} =   a_{i,j+1}, ~~  \sigma_{j} b_{ij} \sigma_{j}^{-1} =   b_{i,j+1}~~~\mbox{for}~1 \leq i < j \leq n-1. \label{ct61}
\end{align}

\bigskip

\section{Subgroups of camomile type} \label{cam}

During last decades people study groups which look like  the braid group $B_n$ (see \cite{B, BB, BBD, BJW, Bir, Ka, KL, V3}).
For example, 

-- Virtual braid group $VB_n$,

-- Welded  braid group $WB_n$,

-- Flat virtual braid group $FVB_n$,

-- Unrestricted virtual braid group $UVB_n$,

-- Singular  braid group $SG_n$,

-- Universal   braid group $UB_n$.

In this list $VB_n$ is generated by $B_n$ and $S_n$. There are epimorphisms
$$
VB_n \to WB_n \to FVB_n \to UVB_n.
$$
On the other side, for $UB_n$ we have two epimorphisms
$$
UB_n \to VB_n,~~UB_n \to SG_n.
$$
Any group from the list above has a pure subgroup, which is the kernel of an epimorphism onto $S_n$. As rule the structures of these pure subgroups are different for the case $n=3$ and $n>3$. In the second case there are far commutativity relations which we do not have in the first case. Furthermore, set of generators and defining relations any pure group for  $n>4$ can be gotten  from the set of generators and defining relations of the pure group for  $n=4$, using conjugations by elements of the corresponding braid group. Hence, in some sense the case $n=4$ is crucial. 
In this section we formalize this observation and suggest a general construction.

Suppose that a group $G$ is defined by a set of generators $X$ and a set of defining relations $R$. In this case we say that $G$ has a presentation
$$
\mathcal{P}(G) = \langle X~|~R \rangle.
$$
We will denote the set of generators $X$ by $\mathcal{G}(G)$ and the set of defining relations $R$ by $\mathcal{R}(G)$.

\begin{definition}
Let $G$ be a group,  $H$  its normal subgroup, $M$ a set of  coset representatives $H$ in $G$. The subgroup  $H$  is said to be a subgroup of {\it camomile type} with  {\it highlighted petal} $H_0$, and a {\it conjugated set} $M_0$,  if   $H_0 \lneqq H$, $M_0$ is a subset of  $M$ that contains the unit element $e$ of $G$, and
$$
\mathcal{P}(H) = \bigcup_{m \in M_0} \mathcal{P}(H_0^m).
$$
It means that
$$
\mathcal{G}(H) = \bigcup_{m \in M_0} \mathcal{G}(H_0^m)~~\mbox{and}~~\mathcal{R}(H) = \bigcup_{m \in M_0} \mathcal{R}(H_0^m).
$$
We will write $H = Cam_G(H_0, M_0)$ and call the subgroups $H_0^m$ by petals.
\end{definition}

\begin{remark}
To avoid trivial cases, we suggest in this definition that the highlighted petal $H_0$ is not equal to $H$.
\end{remark}

The next proposition shows that the highlighted petal is not unique.

\begin{proposition} \label{dif}
Let $G_n = F_n \rtimes S_n$, $n \geq 2$,  be the semi-direct product of the  free group $F_n = \langle x_1, x_2, \ldots, x_n  \rangle$ and the symmetric group $S_n$, which acts on the generators of $F_n$ by the rule $x_i^s = x_{s(i)}$ for $s \in S_n$. Then

1) $F_n$ is a subgroup of   camomile type in $G_n$ with a highlighted petal $F_1 = \langle x_1 \rangle$ and a conjugated set $M_0 = \{ e, (12), (13), \ldots, (1n) \} \subset S_n$;

2) $F_n$ is a subgroup of  camomile type in $G_n$ with a highlighted petal $F_{n-1} = \langle x_1, x_2, \ldots, x_{n-1} \rangle$ and a conjugated set   $M_0 = \{ e, (n-1,n) \}$.
\end{proposition}

\begin{proof}
As the set of coset representatives $F_n$ in $G_n$ we can take $S_n$.

1) Put $M_0 = \{ e, (12), (13), \ldots, (1n) \} \subset S_n$. Since $x_i = x_1^{(1i)}$, $i = 2, 3, \ldots, n$, we have
$$
\mathcal{G}(F_n)  = \bigcup_{m \in M_0} \mathcal{G}(F_1^m) = \{ x_1, x_2, \ldots, x_{n} \}.
$$
Since free group does not contains non-trivial relations,
$$
\mathcal{P}(F_n)  = \bigcup_{m \in M_0} \mathcal{P}(F_1^m).
$$

2) Put $M_0 = \{ e, (n-1,n) \}$. Then $\mathcal{G}(F_{n-1}) = \{ x_1, x_2, \ldots, x_{n-1} \}$ and $\mathcal{G}(F_{n-1}^{(n-1,n)}) = \{ x_1, x_2, \ldots, x_{n-2}, x_n \}$.
Hence,
$$
\mathcal{G}(F_n)  = \mathcal{G}(F_{n-1}) \cup \mathcal{G}(F_{n-1}^{(n-1,n)}).
$$

\end{proof}

Some  groups which are similar to the braid group contain the pure subgroup as subgroup of camomile type.

\begin{example}
1) The pure braid group $P_n$, $n \geq 5$ is a subgroup of camomile type, $P_n = Cam_{B_n}(P_4, M_{n,4})$ with the  highlighted  petal $P_4$ and conjugated set $M_{n,4}$. In this case
$$
\mathcal{P}(P_n) =  \bigcup_{m \in M_{n,4}} \mathcal{P}(P_4^m).
$$
For more details see \cite{BW}.

2) The virtual braid group $VB_n$ is generated by $B_n$ and $S_n$. The virtual pure braid group $VP_n$, $n \geq 5$, is a subgroup of camomile type, $VP_n = Cam_{VB_n}(VP_4, M_{n,4})$ with the   highlighted petal $VP_4$. 
For more details see \cite{BW}.

3) Since there exist epimorphisms 
$$
VP_n \to FVP_n \to UVP_n,
$$
the subgroup $FVP_n \leq FVB_n$ and $UVP_n \leq UVB_n$ are subgroups of camomile type for all $n \geq 5$.
\end{example}

\begin{remark}
A subgroup $H \leq G$ can be a subgroup of camomile type for non-isomorphic   highlighted petal subgroups. For example, $VP_n$, $n \geq 5$, is a subgroup of camomile type with the   highlighted petal $VP_k$ for any $4 \leq k \leq n-1$.
\end{remark}

For singular braid groups the following theorem is true.

\begin{theorem}
The  singular pure braid group $SP_n$, $n \geq 5$, is a group of camomile type, $SP_n = Cam_{SG_n}(SP_4, M_{n,4})$ with the   highlighted petal $SP_4$.
\end{theorem}

\begin{proof}
Using induction by $n$. Let us show that for $n > 5$ the group $SP_n$ is a group of camomile type with the   highlighted petal $SP_{n-1}$ and conjugated set $M_{n,n-1}$. By induction from this it will be follows the need result, i.e.
$$
\mathcal{P}(SP_{n}) = \bigcup_{m \in M_{n,4}} \mathcal{P}(SP_4^m).
$$
We know that
$$
SP_{n-1} = \langle a_{ij}, b_{ij}~|~1 \leq i < j \leq n-1 \rangle,
$$
$$
SP_{n} = \langle a_{ij}, b_{ij}~|~1 \leq i < j \leq n \rangle,
$$
and
$$
M_{n,n-1} =  \{ e, \,  \sigma_{n-1}^{-1},  \, \sigma_{n-1}^{-1} \sigma_{n-2}^{-1},  \,  \ldots,  \, \sigma_{n-1}^{-1} \sigma_{n-2}^{-1} \ldots \sigma_1^{-1} \}.
$$
Since,
$$
 \mathcal{G}(SP_n) =  \mathcal{G}(SP_{n-1}) \cup \{ a_{1n}, b_{1n}, a_{2n}, b_{2n}, \ldots, a_{n-1,n}, b_{n-1,n} \},
$$
it is need to get generators, which do not lie in $SP_{n-1}$. By (\ref{ct61})
$$
 \sigma_{n-1} a_{i,n-1} \sigma_{n-1}^{-1} =   a_{i,n}, ~~  \sigma_{n-1} b_{i,n-1} \sigma_{n-1}^{-1} =   b_{i,n}~~~\mbox{for}~1 \leq i \leq n-1,
$$
and by (\ref{ct41})
$$
 \sigma_{n-2} \sigma_{n-1} a_{n-2, n-1} \sigma_{n-1}^{-1} \sigma_{n-2}^{-1} =   a_{n-1,n}, ~~   \sigma_{n-2} \sigma_{n-1} b_{n-2, n-1} \sigma_{n-1}^{-1} \sigma_{n-2}^{-1}=   b_{n-1,n}.
$$
Hence,
$$
 \mathcal{G}(SP_n) =  \mathcal{G}(SP_{n-1}) \cup  \mathcal{G}(SP_{n-1}^{\sigma_{n-1}^{-1}}) \cup  \mathcal{G}(SP_{n-1}^{\sigma_{n-1}^{-1} \sigma_{n-2}^{-1} }).
$$

Further, the set of defining relations of $SP_n$ is the union of  the set of defining relations of $SP_{n-1}$ and the set of defining relations which contains generators 
$a_{in}$ or $b_{in}$, 
$$
\mathcal{R}(SP_n) = \mathcal{R}(SP_{n-1}) \cup \left( \mathcal{R}(SP_n) \setminus  \mathcal{R}(SP_{n-1}) \right).
$$
Using the conjugation rules (\ref{ct1})--(\ref{ct61}) one can show that any relation in $\mathcal{R}(SP_n) \setminus  \mathcal{R}(SP_{n-1})$ follows from relations 
of the groups
$$
SP_{n-1}^m,~~~m \in M_{n,n-1}.
$$
This completes the proof. 
\end{proof}

\bigskip

\section{The center of $SP_n$} \label{cent}

It is well-known \cite{Mar} that the center $Z(B_n) = Z(P_n)$ is the infinite cyclic group which is generated by element
$$
\Delta_n = (\sigma_1 \sigma_2 \ldots \sigma_{n-1})^n = a_{12} (a_{13} a_{23}) \ldots (a_{1n} a_{2n} \ldots a_{n-1,n}).
$$

It was shown \cite{FRZ}  (see also \cite{V4}), that  $ Z(SG_n) \cong Z(SP_n) \cong Z(B_n)$. Neshchadim \cite{N1, N2} proved that $Z(P_n)$ is a direct factor in $P_n$.
In \cite{BK} was proved  that $Z(SP_3)$ is a direct factor in $SP_3$. In this section we prove the same result for arbitrary $n \geq 3$. We will use the following notations
$$
\delta_k = a_{1k} a_{2k} \ldots a_{k-1,k},~~~k = 2, 3, \ldots, n.
$$
Then
$$
\Delta_n = \delta_2 \delta_3 \ldots \delta_n.
$$

Using defining relations of $SP_n$ one can prove the next lemma.

\begin{lemma} \label{cl1}
Let $n \geq 3$, then for any $2< k < j \leq n$ the following formulas hold
 \begin{align}
& a_{1j}^{a_{2j}^{-1} a_{1j}^{-1} \delta_3 \delta_{4} \ldots \delta_j} =  a_{1j}, ~~b_{1j}^{a_{2j}^{-1} a_{1j}^{-1} \delta_3 \delta_{4} \ldots \delta_j} =  b_{1j}, \label{cent1}\\
& a_{2j}^{ a_{1j}^{-1} \delta_3 \delta_{4} \ldots \delta_j} =  a_{2j}, ~~b_{2j}^{a_{1j}^{-1} \delta_3 \delta_{4} \ldots \delta_j} =  b_{2j}, \label{cent2},\\
& a_{kj}^{\delta_k \delta_{k+1} \ldots \delta_j} =  a_{kj}, ~~b_{kj}^{\delta_k \delta_{k+1} \ldots \delta_j} =  b_{kj}, \label{cent3}
\end{align}
\end{lemma}

\begin{lemma} \label{cl2}
Let $n \geq 3$, then for any $1< k < j \leq n$ the following formulas hold
$$
(a_{1j} a_{2j} \ldots a_{j-1,j})^{a_{ik}} = a_{1j} a_{2j} \ldots a_{j-1,j},
$$
$$
(a_{1j} a_{2j} \ldots a_{j-1,j})^{b_{ik}} = a_{1j} a_{2j} \ldots a_{j-1,j}.
$$
\end{lemma}

\begin{proof}
The firs  formula holds in $P_n$ (see \cite{Mar}).

Let us prove the second formula. Suppose that $i=k-1$. Then
$$
(a_{1j} a_{2j} \ldots a_{k-1,j} a_{kj} \ldots a_{j-1,j})^{b_{k-1,k}} = (a_{1j}  \ldots a_{k-2,j})^{b_{k-1,k}} (a_{k-1,j} a_{kj})^{b_{k-1,k}}
(a_{k+1, j}  \ldots a_{j-1,j})^{b_{k-1,k}}.
$$
Since generator $b_{k-1,k}$ commutes with all generators
$$
a_{1j},   \ldots,  a_{k-2,j}, a_{k+1, j},   \ldots,  a_{j-1,j},
$$
and from relation (\ref{dr-6}) follows that
$$
(a_{k-1,j} a_{kj})^{b_{k-1,k}} = a_{k-1,j} a_{kj},
$$
we get the need formula.

Suppose that $1 \leq  i < k-1$. Then
$$
\left((a_{1j} \ldots a_{i-1,j}) (a_{ij} \ldots a_{kj}) (a_{k+1,j} \ldots a_{j-1,j})\right)^{b_{ik}} =
$$
$$
~~~~~~~~~~~~~~~~~~~~~~= (a_{1j} \ldots a_{i-1,j})^{b_{ik}} (a_{ij} \ldots a_{kj})^{b_{ik}}
(a_{k+1, j}  \ldots a_{j-1,j})^{b_{ik}}.
$$
As in the previous case, the generator  $b_{ik}$ commutes with all generators
$$
a_{1j},   \ldots,  a_{i-1,j}, a_{k+1, j},   \ldots,  a_{j-1,j}.
$$
Further,
$$
(a_{ij} \ldots a_{kj})^{b_{ik}} = (a_{ij} a_{kj})^{b_{ik}}  (a_{kj}^{-1} a_{i+1,j} a_{kj})^{b_{ik}}  (a_{kj}^{-1} a_{i+2,j} a_{kj})^{b_{ik}}  \ldots  (a_{kj}^{-1} a_{k-1,j}a_{kj})^{b_{ik}}.
$$
Using the relations (\ref{dr-6}) and (\ref{plast}) we get the need formula.

\end{proof}

Now we are ready to prove

\begin{theorem}
For any $n \geq 3$ the center $Z(SG_n)$ is a direct factor in $SP_n$. But $Z(SG_n)$ is not a direct factor in $SG_n$.
\end{theorem}

\begin{proof}
As we know $SG_n$ is generated by
$$
a_{ij}, ~~b_{ij},~~1 \leq i < j \leq n,
$$
and is defined by the set of relations $R$, i.e.
\begin{equation} \label{pres0}
SP_n = \langle a_{ij}, ~~b_{ij},~~1 \leq i < j \leq n,~|~R \rangle.
\end{equation}
The set of relations $R$ is disjoint union of two subsets, $R = R_1 \sqcup R_2$, where $R_1$ is the set of relations which  contain $a_{12}$ and $R_2$ is the set of relations which do not contain  $a_{12}$. Denote by $A$ the set of generators $SP_n$ without generator $a_{12}$ and denote by $H_n = \langle A \rangle \leq SP_n$.

Let us prove  that $SP_n$ also has the following presentation
\begin{equation} \label{pres}
SP_n = \langle A, \Delta_n~|~R_2, [\Delta_n, a] = 1, a \in A \rangle.
\end{equation}
It is enough to prove that any relation from $R_1$ follows from relations $R_2$ and relations  $[\Delta_n, a] = 1$, $a \in A$. Any relation from  $R_1$ has one of the forms:\\

1) $a_{1j}^{a_{12}} = a_{1j}^{a_{2j}^{-1} a_{1j}^{-1}}$  \, \, \,   $b_{1j}^{a_{12}} = b_{1j}^{a_{2j}^{-1} a_{1j}^{-1}}$;\\

2) $a_{2j}^{a_{12}} = a_{2j}^{a_{1j}^{-1}}$   \, \, \,    $b_{2j}^{a_{12}} = b_{2j}^{a_{1j}^{-1}}$;\\

3) $a_{kj}^{a_{12}} = a_{kj}$ \, \, \,   $a_{kj}^{a_{12}} = a_{kj}$  for $2 < k \leq n$.\\

Since $\Delta_n = \delta_2 \delta_3 \ldots \delta_n$ and $\delta_2 = a_{12}$, then $a_{12} = \Delta_n  \delta_n^{-1} \delta_{n-1}^{-1} \ldots \delta_3^{-1}$. Using this formula, we can remove $a_{12}$ from the generating set of $SP_n$. Hence, $SP_n$ is generated by $A$ and $\Delta_n$.

Let us show that we can remove the set of relations $R_1$ and insert the relations  $[\Delta_n, a] = 1, a \in A$. The first relation of the type 1) can be written in  the form
$$
a_{1j}^{\Delta_n  \delta_n^{-1} \delta_{n-1}^{-1} \ldots \delta_3^{-1}} = a_{1j}^{a_{2j}^{-1} a_{1j}^{-1}}.
$$
The element $\Delta_n$ lies in the center of $SP_n$, hence
$$
a_{1j}^{\delta_n^{-1} \delta_{n-1}^{-1} \ldots \delta_{j+1}^{-1}} = a_{1j}^{a_{2j}^{-1} a_{1j}^{-1} \delta_{3} \ldots \delta_j}.
$$
From the first relation of Lemma \ref{cl1} follows that the left side of this relation is equal to $a_{1j}$. From the first relation of (\ref{cent1}) follows that the right side of this relation is equal to $a_{1j}$.
Hence, we can remove the firs relation of the form 1). By the same way, using the second  relation of Lemma \ref{cl1} and the second relation of (\ref{cent1}) we can show
that we can remove the second relation of the form 1).

Analogously, consider the relations of the form 2) and 3) and using Lemmas \ref{cl1} and \ref{cl2} we can prove that $SP_n$ has the presentation (\ref{pres}).
From this  presentation  follows that there are two epimorphisms
$$
\pi_1 : SP_n \to Z(SG_n),~~\pi_1(\Delta_n) = \Delta_n,~~\pi_1(a) = 1 ~\mbox{for all}~a \in A;
$$
$$
\pi_2 : SP_n \to H_n,~~\pi_2(\Delta_n) = 1,~~\pi_2(a) = a ~\mbox{for all}~a \in A.
$$
Hence, $SP_n = \langle Z(SG_n), H_n \rangle$, the subgroup $H_n$ has a presentation
$$
H_n = \langle A~|~R_2 \rangle
$$
and $Z(SG_n) \cap H_n = 1$. We proved   the first part of the theorem.

The second part of the theorem follows from the fact that there exists an epimorphism $SG_n \to B_n$ and from the fact that $Z(B_n)$ is not a direct factor of $B_n$.
\end{proof}

\bigskip

\section{Singquandles and invariant  of singular links}\label{Gr}


A singular link in $\mathbb{S}^3$ is the image of a smooth immersion of $n$ circles in $\mathbb{S}^3$ that has finitely many double points, called singular points. An orientation of each circle induces orientation on each component of the link. This gives an oriented singular link. In the present section we consider only oriented singular links and oriented singular link diagrams. A singular link diagram is a  projections of the singular link to the plane such that any crossing is a classical crossing with over and under information, or a singular crossing.  Two oriented singular link diagrams are  equivalent if and only if one can obtain  from the
other by a finite sequence of singular Reidemeister moves (see, for example, \cite{E}).

In the  classical knot theory for any link diagram it can be defined a quandle which is a link invariant (see \cite{M}, \cite{J}). 
 Recall that an algebraic system  is a non-empty set with some algebraic operations.
A quandle is an algebraic system $Q$ with a binary operation $(x,y) \mapsto x * y$ satisfying the following axioms:
\begin{enumerate}
\item[(Q1)] $x*x=x$ for all $x \in Q$,
\item[(Q2)] For any $x,y \in Q$ there exists a unique $z \in Q$ such that $x=z*y$,
\item[(Q3)] $(x*y)*z=(x*z) * (y*z)$ for all $x,y,z \in Q$.
\end{enumerate}
\begin{figure}[ht]
\includegraphics[width=12cm, height=5cm]{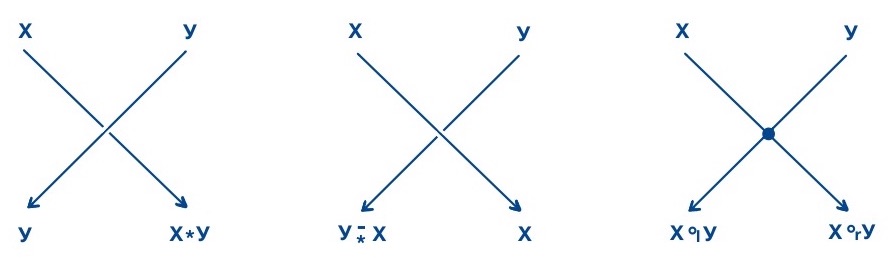}
\caption{Classical and singular crossings} \label{figure2}
\end{figure}

\begin{figure}[h]
\includegraphics[width=5.4cm, height=5.4cm]{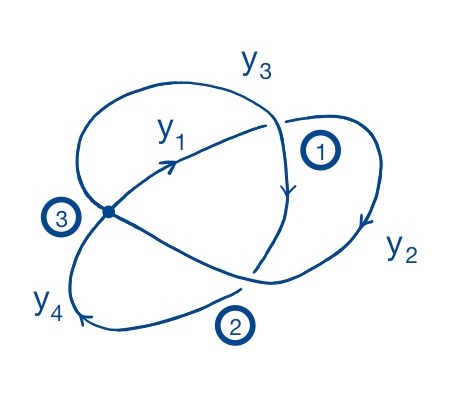}\\
\includegraphics[width=16cm, height=5.5cm]{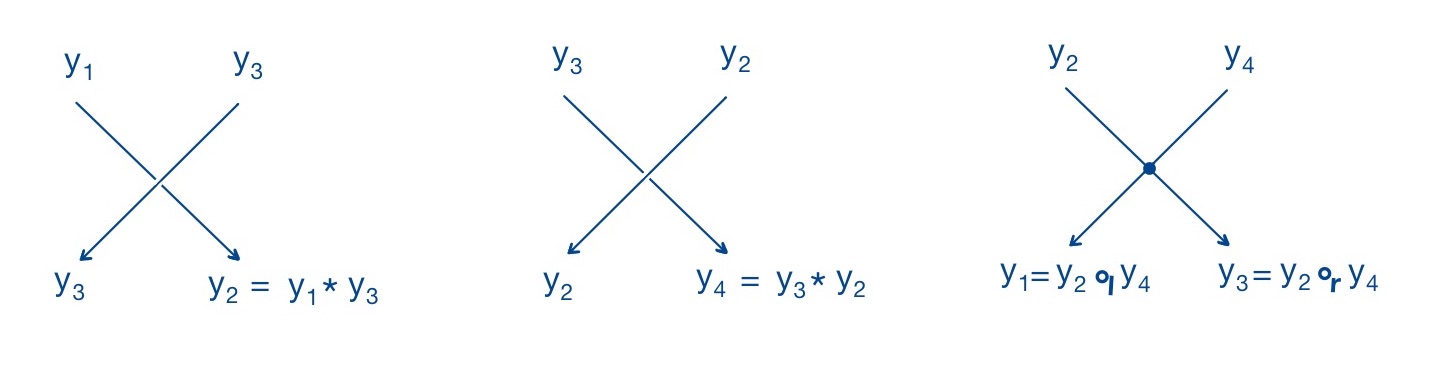}\\
\caption{Singular trefoil and relations} \label{figure3}
\end{figure}

By analogy with the definition of quandle, in \cite{E} was introduced  an {\it oriented singquandle} whose axioms  come  from singular  Reidemeister moves on oriented singular link diagram. We give a definition of oriented singquandle in the form that is different
from  \cite{E}, but it is easy to see that our definition is equivalent to the original one.

\begin{definition}
An oriented singquandle $(X, *, \circ_l, \circ_r)$ is an algebraic system with three binary operations  such that $(X, *)$ is a quandle and for any $x, y, z \in X$ hold
 \begin{align}
& (x  \circ_l y) * z = (x * z) \circ_l (y * z),   \label{re1}\\
& (x  \circ_r y) * z = (x * z) \circ_r (y * z),   \label{re1}\\
& (x  \circ_l y) \bar{*} z = (x \bar{*} z) \circ_l (y \bar{*} z),   \label{re2}\\
& (x  \circ_r y) \bar{*} z = (x \bar{*} z) \circ_r (y \bar{*} z),   \label{re2}\\
& (y \bar{*} (x \circ_l z)) * x = (y * (x \circ_r z))  \bar{*} z,   \label{re3}\\
& x \circ_r y = y \circ_l ( x * y),   \label{re4}\\
& (x \circ_l y) * (x \circ_r y)  = y \circ_r (x * y).   \label{re5}
\end{align}
\end{definition}

\begin{remark}
If we put  $z = x * y$, then  $x = z \bar{*} y$ and (\ref{re4}) has the form
$$
y \circ_l z = (z  \bar{*} y) \circ_r y.
$$
Hence, in the definition of singquandle we can remove one operation and present it as a quandle with one extra operation (see \cite{BC}).
\end{remark}


Now we recall a definition of fundamental oriented singquandle, following  \cite{E}. Let $D_L$ be a diagram of an  oriented singular link $L$. Removing singular points of $D_L$ we get some number of connected component which we call by arcs of the diagram. Letting these arcs by $y_1, y_2, \ldots, y_m$. Then the fundamental
oriented singquandle  $SQ(D_L)$ is generated by $y_1, y_2, \ldots, y_m$ and defined by the set of relations from the Figure \ref{figure2} for any crossing (classical and singular). As was proved in \cite{E} this singquandle $SQ(D_L)$ is an invariant of $L$ and we will denote it as $SQ(L)$. If $L$ is a classical link, then $SQ(L)$ is the fundamental quandle of $L$ that is defined by Matveev \cite{M} and independently Joice \cite{J}. In the case of classical knots the fundamental quandle is  almost complete invariant, it means that if $K$ and $K'$ are two classical knots and their quandles are isomorphic, then $K$ is equivalent to $K'$ or to $-\overline{K'}$, where minus means the knot with opposite orientation and bar means the mirror image.

\begin{example} \label{ex1}
Let $T_s$ be  the singular trefoil knot. Its diagram is depictured in  Figure  \ref{figure3}.  After removing the singular crossing we get four arcs, which we letting by $y_1$, $y_2$, $y_3$, and $y_4$. These are generators of the singquandle $SQ(T_s)$. In the first crossing of the diagram we have relation $y_2 = y_1 * y_3$, is the second crossing we have relation
$y_4 = y_3 * y_2$ and in the third crossing, which is the singular crossing we have two relations,
$$
y_1 = y_2 \circ_l y_4,~~~y_3 = y_2 \circ_r y_4.
$$
Hence,
$$
SQ(T_s) = \langle y_1, y_2, y_3, y_4 ~|~  y_2 = y_1 * y_3,~~~ y_4 = y_3 * y_2~~~ y_1 = y_2 \circ_l y_4,~~~y_3 = y_2 \circ_r y_4 \rangle.
$$
We can remove the generator $y_2$ and $y_4$, then
$$
SQ(T_s) = \langle y_1,  y_3 ~|~  y_1 =(y_1 * y_3) \circ_l (y_3 * y_2),~~~y_3 = (y_1 * y_3) \circ_r (y_3 * y_2) \rangle.
$$
\end{example}

Now we give other  interpretation of $SQ(L)$, using a representation of the singular braid monoid $SB_n$. Here we are following the ideas from \cite{BMN} and \cite{BMN-1},
where were defined groups of virtual links.

Let $FSQ_n$ be a free oriented singquandle that is generated by $x_1, x_2, \ldots, x_n$. Define a map $\Phi : SB_n \to \End(FSQ_n)$ which is defined on the generators by the rules
$$
\Phi(\sigma_{i}) : \left\{
\begin{array}{ll}
x_{i} \longmapsto  x_{i+1}, &  \\
x_{i+1} \longmapsto x_{i} * x_{i+1}, & \\
x_{j} \longmapsto x_{j},   j\neq i,i+1, &
\end{array} \right.~~~~
\Phi(\tau_{i}) : \left\{
\begin{array}{ll}
x_{i} \longmapsto x_{i} \circ_l x_{i+1}, &  \\
 x_{i+1} \longmapsto x_{i} \circ_r x_{i+1}, & \\
x_{j} \longmapsto x_{j},   j\neq i,i+1. &
\end{array} \right.
$$
Using a standard technique it is easy to prove
 
\begin{proposition}
The map $\Phi : SB_n \to \End(FSQ_n)$ is a representation.
\end{proposition}

Hence, using induction on the length, we can define $\Phi$ on arbitrary word $w$ which presents some element of $SB_n$. If $w = u \diamond  v$, where $\diamond  \in \{ *, \bar{*}, \circ_l, \circ_r \}$ and the length of $u$ and $v$ less of the length of $w$, then we define $\Phi(w) = \Phi(u) \diamond \Phi(u)$.

\medskip

Since any oriented singular link is the closure of some singular braid \cite{Bir}, then any singular link $L$ is the closure of a singular braid $\beta \in SB_n$. Using the same approach as in \cite{BMN}, one can give another interpretation of fundamental oriented singquandle.

\begin{theorem}
Suppose that a singular link $L$ is the closure of the singular braid $\beta \in SB_n$. Then
$$
SQ(L) \cong \langle x_1, x_2, \ldots, x_n~|~x_i = \Phi(\beta)(x_i),~~i = 1, 2, \ldots, n \rangle,
$$
where we consider the action of $\Phi(\beta)$ on $FSQ_n$ from left to right.
\end{theorem}

\begin{example}
Let $\beta = \sigma_1^2 \tau_1 \in SB_3$. We find the image $\Phi( \sigma_1^2 \tau_1)$, acting on the generators $x_1$, $x_2$ from the right to the left,
$$
\Phi(\sigma_1^2 \tau_1) : \left\{
\begin{array}{ll}
x_{1} \mapsto  (x_{1} * x_2) \circ_l (x_2 * (x_1 * x_2)), &  \\
x_{2} \mapsto  (x_{1} * x_2) \circ_r (x_2 * (x_1 * x_2)). &
\end{array} \right.
$$
Define a sigquandle
$$
SQ_{\Phi}(\sigma_1^2 \tau_1) = \langle x_1, x_2~|~x_1 =  (x_{1} * x_2) \circ_l (x_2 * (x_1 * x_2)),~~x_2 = (x_{1} * x_2) \circ_r (x_2 * (x_1 * x_2)) \rangle.
$$
One can see that it is the singquandle of the singular trefoil knot which we constructed in Example \ref{ex1}.
\end{example}

\medskip

\section{Representation of $SB_n$ by endomorphisms of free group and groups of singular links} \label{GSL}

In this section we construct some groups of singular links which are link invariants. 
To do it we use the next representations (see \cite{E})  of $SB_n$ into $\End(F_n)$, which are extensions of the Artin representation $\varphi_A \colon B_n \to \Aut(F_n)$ (see Section \ref{BD}):  

$$\Phi_1({\tau _i}):\left\{ {\begin{array}{*{20}{c}}
{{x_i} \to {x_i}{x_{i+1}}{x_i}x_{i+1}^{ - 1}x_i^{ - 1},}\\
{{x_{i+1}} \to {x_i}{x_{i+1}}x_i^{ - 1},\;\;\;\;\;\;\;}\\
{{x_j} \to {x_j},~~~~j \not= i, i+1,\;\;\;\;\;}
\end{array}} \right.\;\;\;\;\;\;\;\Phi_1({\sigma _1}):\left\{ {\begin{array}{*{20}{c}}
{{x_i} \to {x_i}{x_{i+1}} x_i^{ - 1},\;\;\;\;\;\;\;}\\
{{x_{i+1}} \to {x_i},\;\;\;\;\;\;\;\;\;\;\;\;\;\;\;}\\
{{x_j} \to {x_j},~~~~j \not= i, i+1,\;\;\;}
\end{array}} \right.$$

$$
\Phi_2(\tau_{i}) : \left\{
\begin{array}{ll}
x_{i} \longmapsto x_{i+1}^{-1} \, x_{i} \, x_{i+1}, &  \\
 x_{i+1} \longmapsto x_{i+1}^{-1} x_{i}^{-1}x_{i+1} x_i x_{i+1}, & \\
  x_{j} \longmapsto x_{j},   j \neq i,i+1, &
\end{array} \right.~~~~
\Phi_2(\sigma_{i}) : \left\{
\begin{array}{ll}
x_{i} \longmapsto x_{i} \, x_{i+1} \, x_i^{-1}, &  \\
x_{i+1} \longmapsto x_{i}, & \\
x_{j} \longmapsto x_{j},   j\neq i,i+1, &
\end{array} \right.
$$

$$
\Phi_3(\tau_{i}) : \left\{
\begin{array}{ll}
x_{i} \longmapsto x_i x_{i+1}^{-1} \, x_{i}^{-1} \, x_{i+1} \, x_i, &  \\
 x_{i+1} \longmapsto x_{i}^{-1} x_{i+1}^{-1} \, x_{i}  \, x_{i+1}^2, & \\
x_{j} \longmapsto x_{j},   j\neq i,i+1, &
\end{array} \right.~~~~
\Phi_3(\sigma_{i}) : \left\{
\begin{array}{ll}
x_{i} \longmapsto x_{i} \, x_{i+1} \, x_i^{-1}, &  \\
x_{i+1} \longmapsto x_{i}, & \\
x_{j} \longmapsto x_{j},   j\neq i,i+1, &
\end{array} \right.
$$

$$
\Phi_{4,n}(\tau_{i}) : \left\{
\begin{array}{ll}
x_{i} \longmapsto x_{i+1}^{-1} \, (x_{i}^{-1} \, x_{i+1})^n, &  \\
 x_{i+1} \longmapsto (x_{i+1}^{-1} x_{i})^{n+1} \, x_{i+1}, & \\
x_{j} \longmapsto x_{j},   j\neq i,i+1, &
\end{array} \right.~~~~
\Phi_{4,n}(\sigma_{i}) : \left\{
\begin{array}{ll}
x_{i} \longmapsto x_{i} \, x_{i+1} \, x_i^{-1}, &  \\
x_{i+1} \longmapsto x_{i}, & \\
x_{j} \longmapsto x_{j},   j\neq i,i+1, &
\end{array} \right.
$$
where $n \geq 1$.

\begin{remark}
 In the paper \cite{BKV} were constructed some representations of the subgroup 
$T_n = \langle \tau_1, \tau_2, \ldots, \tau_{n-1} \rangle \leq SG_n$  into $\Aut(F_n)$ and proved that some of them are faithful.
\end{remark}

Using the Artin representation $\varphi_A$ it is possible to define a group $G_A(\beta)$, $\beta \in B_n$ and prove that this group is an invariant of the link $L=\hat{\beta}$ that is the closure of $\beta$. In fact, this group is the fundamental group of the compliment of $L$ in the 3-dimensional space, $G_A(\beta) = \pi_1(\mathbb{R}^3 \setminus L)$ (see \cite{Bir1}).

Let us define some set of groups for a singular braid.

\begin{definition}
Let $\beta \in SB_n$ and $\Psi : SB_n \to \End(F_n)$ be a representation from the list $\Phi_1, \Phi_2, \Phi_3, \Phi_{4,n}$, $n \geq 1$, then we define a group
$$
G_{\Psi}(\beta) = \langle x_1, x_2, \ldots, x_n~|~x_i = \Psi(\beta)(x_i),~~i = 1, 2, \ldots, n \rangle.
$$
For simplicity we shall  write $G_{1}(\beta)$, $G_{2}(\beta)$, $G_{3}(\beta)$, $G_{4,n}(\beta)$ instead  $G_{\Phi_1}(\beta)$, $G_{\Phi_2}(\beta)$, $G_{\Phi_3}(\beta)$, $G_{\Phi_{4,n}}(\beta)$,  correspondingly.
\end{definition}

A proof of the next theorem is the same as for  virtual links (see, for example \cite{BMN, BMN-1}).

\begin{theorem}
Suppose that a singular link $L$ is the closure of the singular braid $\beta \in SB_n$. Then any group
$G_{\Psi}(\beta)$
is an invariant of $L$ that means that it does not depends on the realization of $L$ by a singular braid.
\end{theorem}

Since for any classical braid $\beta$ we have
$$
\Phi_1(\beta) =  \Phi_2(\beta) = \Phi_3(\beta) = \Phi_{4,n}(\beta) = \varphi_A(\beta), 
$$
then
$$
G_{1}(\beta) = G_{2}(\beta) = G_{3}(\beta) = G_{4,n}(\beta) = G_A(\beta)
$$
is the classical group of the link $L = \hat{\beta}$.
For example, any of these groups for  
 trivial $m$-component link is the free group of rank $m$. 
The Hopf link can be presented as  the closure of the braid $\sigma_1^2 \in B_2$ and its group is the free abelian group of rank 2.

The singular Hoph link is the closure of the singular braid $\sigma_1 \tau_1$. Using the representations $\Phi_1$, $\Phi_2$, $\Phi_3$, $\Phi_{4,n}$, we will find the groups
$G_1$, $G_2$, $G_3$, $G_{4,n}$, for the singular Hopf link $L_{sf}$ and for its mirror image $\overline{L_{sf}}$ that is the closure of the braid $\sigma_1^{-1} \tau_1$.

1) For $G_1$ we  have
$$
\Phi_1(\sigma_1 \tau_1 ) : \left\{
\begin{array}{ll}
x_{1} \longmapsto  x_{1} x_2 x_1 x_2 x_1^{-1}  x_2^{-1}  x_1^{-1}, &  \\
x_{2} \longmapsto x_{1} x_2 x_1  x_2^{-1}  x_1^{-1}. &
\end{array} \right.
$$
Then
$$
G_1(L_{sf}) = \langle x_1, x_2~|~x_1 = x_2 x_1  x_2 x_1^{-1}  x_2^{-1},~~x_2 = x_1 x_2  x_1 x_2^{-1}  x_1^{-1} \rangle.
$$
This group is the braid group $B_3$ that is the group of the trefoil knot.

For the mirror image  $G_1(\overline{L_{sf}})$ is the infinite cyclic group.

2) Using the similar calculations it is easy to check that
$$
G_2(L_{sf}) \cong \mathbb{Z},~~~G_2(\overline{L_{sf}}) \cong B_3.
$$

3) For $G_3$ we have
$$
\Phi_3(\sigma_1 \tau_1 ) : \left\{
\begin{array}{ll}
x_{1} \longmapsto  x_{1} x_2  x_1^{-1}  x_2^{-1} x_1 x_2  x_1^{-1}, &  \\
x_{2} \longmapsto  x_1  x_2^{-1}  x_1^{-1} x_2 x_1. &
\end{array} \right.
$$
Then
$$
G_3(L_{sf}) = \langle x_1, x_2~|~x_1 = x_2  x_1^{-1}  x_2^{-1} x_1  x_2,~~x_2 = x_1  x_2^{-1}  x_1^{-1} x_2 x_1 \rangle \cong
\langle x_1, x_2~|~x_1   x_2^{-1}  x_1^{-1} =  x_2   x_1^{-1} x_2^{-1} \rangle
$$
and
$$
G_3(\overline{L_{sf}}) = \langle x_1, x_2~|~x_2 x_1^2 =   x_1 x_2^2 \rangle.
$$

4) For $G_{4,n}$ we  have
$$
\Phi_{4,n}(\sigma_1 \tau_1 ) : \left\{
\begin{array}{ll}
x_{1} \longmapsto   x_2^{-1} ( x_1^{-1} x_2)^n (x_2^{-1} x_1)^{n+1} x_2 (x_2^{-1} x_1)^n x_2, &  \\
x_{2} \longmapsto   x_2^{-1} ( x_1^{-1} x_2)^n. &
\end{array} \right.
$$
Then
$$
G_{4,n}(L_{sf}) = \langle x_1, x_2~|~x_1 = x_2^{-2}  x_1 x_2 (x_2^{-1}  x_1)^n  x_2,~~x_2^2 = (x_1^{-1} x_2)^n \rangle \cong
$$
$$
\langle x_1, x_2~|~x_1 = x_2^{-2}  x_1^2 (x_2^{-1}  x_1)^{n-2}  x_2,~~x_2 = (x_1^{-1} x_2)^{n-1} x_1^{-1} \rangle.
$$
In particular, $G_{4,1}(L_{sf}) \cong \mathbb{Z}_2$.

For the mirror image,
$$
G_{4,n}(\overline{L_{sf}}) = \langle x_1, x_2~|~x_1 = (x_2^{-1}  x_1)^{n-1} x_2^{-1}  x_1  x_2^{-1} x_1 x_2,~~  
x_2 = (x_1^{-1}  x_2)^{n-1} x_1^{-1}  x_2  x_1^{-1} x_2^{-1} x_1\rangle.
$$

\bigskip

\section{Some open problems and directions for further research}

\subsection{How strong invariant is the fundamental singquandle?}

As we know, fundamental quandle is strong invariant. It is natural to formulate

\begin{problem}
Suppose that $K$ and $K'$ are two oriented singular knots such that $SQ(K) \cong SQ(K')$. What can we say on $K$ and $K'$?
\end{problem}

\subsection{Comparing of group $SP_n$ with group $ST_n$}

In \cite{GKM} was defined a group $ST_n = Ker(SG_n \rightarrow  S_n)$, where the homomorphism acts  by the rules:
$$
\sigma_i \mapsto e,~~~\tau_i \mapsto (i, i+1),~~~i = 1, 2, \ldots, n-1. 
$$
It has been proven that  $ST_3$ is isomorphic to  $SP_3$.
\begin{problem}
Is it true that $ST_n$ is isomorphic to the pure singular braid group $SP_n$, for all $n > 3?$
\end{problem}

\end{document}